# Why we must heed Wittgenstein's "notorious paragraph"

Bhupinder Singh Anand

In this paper, we argue that, although Wittgenstein's reservations on Gödel's interpretation of his own formal reasoning are, indeed, of historical importance, the uneasiness that academicians and philosophers continue to sense, and express, over standard interpretations of Gödel's formal reasoning - even seventy years after the publication of his seminal 1931 paper - is of much greater significance, and relevance, to us today.

## Contents

## 1. Introduction

In their paper "A note on Wittgenstein's 'notorious paragraph' about the Gödel Theorem"[1] [FP00], Juliet Floyd and Hilary Putnam draw attention to Wittgenstein's remarks[2] ([Wi78], Appendix III §8):

> I imagine someone asking my advice; he says: "I have constructed a proposition (I will use 'P' to designate it) in Russell's symbolism, and by means of certain definitions and transformations it can be so interpreted that it says: 'P is not provable in Russell's system.' Must I not say that this proposition on the one hand is true, and on the other hand is unprovable? For suppose it were false; then it is true that it is provable. And that surely cannot be! And if it is proved, then it is proved that it is not provable. Thus it can only be true, but unprovable."

---

[1] Downloadable PDF file.

[2] In footnote 9, Floyd and Putnam [FP00] note that: The "notorious" paragraph RFM I Appendix III §8 was penned on 23 September 1937, when Wittgenstein was in Norway (see the Wittgenstein papers, CD Rom, Oxford University Press and the University of Bergen, 1998, Item 118 (Band XIV), pp. 106ff).

> Just as we ask, "'Provable' in what system?," so we must also ask, "'True' in what system?" "True in Russell's system" means, as was said, proved in Russell's system, and "false in Russell's system" means the opposite has been proved in Russell's system.--Now what does your "suppose it is false" mean? In the Russell sense it means, "suppose the opposite is proved in Russell's system"; if that is your assumption you will now presumably give up the interpretation that it is unprovable. And by "this interpretation" I understand the translation into this English sentence. -- If you assume that the proposition is provable in Russell's system, that means it is true in the Russell sense, and the interpretation "P is not provable" again has to be given up. If you assume that the proposition is true in the Russell sense, the same thing follows. Further: if the proposition is supposed to be false in some other than the Russell sense, then it does not contradict this for it to be proved in Russell's system. (What is called "losing" in chess may constitute winning in another game.)

Floyd and Putnam then argue that this paragraph contains a "philosophical claim of great interest which":

> ... is simply this: if one assumes (and, a fortiori if one actually finds out) that ~P is provable in Russell's system one should ... give up the "translation" of P by the English sentence "P is not provable".

In this paper, we argue that, although Wittgenstein's reservations on Gödel's interpretation of his own formal reasoning are, indeed, of historical importance, the uneasiness that academicians and philosophers such as Floyd and Putnam continue to sense, and express, over standard (text-book) interpretations of Gödel's formal reasoning - even seventy years after the publication of his seminal 1931 paper [Go31a] - is of much greater significance, and relevance, to us today.

In a set of related, unpublshed, papers[3], we argue that standard interpretations of foundational concepts of classical mathematical theory may be implicitly influenced by, and built upon, some of Gödel's questionable interpretations of his own formal reasoning[4]; hence such interpretations can be argued as being, prima facie, either ambiguous, or non-constructive, or both.

We now argue further that Wittgenstein's reservations, and Floyd and Putnam's uneasiness, can - and arguably must, as we advocate in this paper - be seen as indicating specific points of such ambiguity that need to be addressed on philosophical grounds, rather than dismissed on technicalities.

## 2. What is mathematics?

### 2.1 Mathematics as a set of formal languages

Without attempting to address the issue in its broader dimensions, we take Wittgenstein's remarks as implicitly indicating that mathematics is to be considered simply as a set of formal languages.

Any language of such a set, say Peano Arithmetic PA[5] (or Russell and Whitehead's Principia Mathematica, PM, or ZFC[6]), expresses - in a finite, unambiguous, and communicable manner - relations between concepts[7] that are external to the language PA (or to PM, or to ZFC). Each such language is, thus, essentially two-valued[8], since a relation[9] either holds or does not hold externally (relative to the language).

---

[3] [An02a] to [An03g].

[4] For instance, see [An02c], §4.4 and §4.5, and §5 of [An03b].

[5] By "PA", or "standard PA", we mean the formal first order theory S defined by Mendelson ([Me64], p102-3).

[6] We assume that "PM" and "ZFC" are familiar terms; the precise definitions of these formal languages is not relevant to the arguments of this paper.

[7] For reasons that should become clearer later, we prefer the term "concepts" to the term "entities"; we intend to define some mathematical concepts as "mathematical objects", where some "concept" may not represent what we would intuitively term as an "entity".

Further, a selected, finite, number of primitive formal assertions[10] about a finite set of selected primitive relations of, say, PA are defined as axiomatically PA-provable; all other assertions about relations that can be effectively[11] defined in terms of the primitive relations are termed as PA-provable if, and only if, there is a finite sequence of assertions of PA, each of which is either a primitive assertion, or which can effectively be determined in a finite number of steps as an immediate consequence of any two assertions preceding it in the sequence by a finite set of rules of consequence.

**2.2 An interpretation must be an effectively decidable translation**

An effective interpretation[12] of a language such as PA into another language, say PM (or ZFC, or even English), is essentially the specification of an effective method by which any assertion of PA is translated unambiguously into a unique assertion of PM (or ZFC, or English). Clearly, if a relation is provable in PA, then it should be effectively decidable in any interpretation of PA - since a finite proof sequence of PA would, prima facie, translate as a finite proof sequence in the interpretation.

---

[8] We discuss the significance of this later in the section §4.2 on omega-consistency.

[9] The term "relation" is thus treated as a primitive, undefined mathematical concept.

[10] We prefer the term "formal assertion" to the familiar terms "formal sentence", or "formal proposition", since a paradigm shift is involved in interpreting the expression "[$(Ax)F(x)$]" constructively. For convenience, we may also refer to a "formal assertion" as an "assertion", where the intention is clear from the context.

[11] "Effectively" means by some finite, unambiguous, mechanical procedure (cf. [Me64], p 207).

[12] We sometimes use the term "interpretation" as a noun in its commonly understood sense, and sometimes in a mathematical sense as in this instance; Mendelson ([Me64], §2, p49) gives the more precise, classical, definition of the term when used in a mathematical sense. Unless specified otherwise, whether the term is to be understood generally, or in its mathematical sense, is determined by the context.

## 2.3 Is the converse necessarily true?

The question arises: If an assertion is decidable in an interpretation M of PA, then does such decidability translate into an effective method of decidability in PA?

Obviously, such a question can only be addressed unambiguously if there is an effective method for determining whether an assertion is decidable in M. If there is no such effective method, then we are faced with the following thesis that is implicit in, and central to, Wittgenstein's "notorious" remark:

> If there is no effective method for the unambiguous decidability of assertions in a language, then it is not a meaningful mathematical language.

In other words, in the absence of an effective method of decidability in M, any correlation of a PA-assertion with an assertion in M is essentially arguable, and it is meaningless to ask whether an assertion of PA is decidable under interpretation in M or not (the question of whether the assertion is decidable in PA or not is, then, an issue of secondary consequence).

## 2.4 Tarskian truth under the standard interpretation

The philosophical dimensions of this thesis emerge if we take M as the standard, arithmetical, interpretation of PA, where [Me64]:

(*a*) the set of non-negative integers is the domain,

(*b*) the integer 0 is the interpretation of the symbol "0" of PA,

(*c*) the successor operation (addition of 1) is the interpretation of the " ' " function (i.e. of $f_1^1$),

(*d*) ordinary addition and multiplication are the interpretations of "+" and ".",

(*e*) the interpretation of the predicate letter "=" is the equality relation.

Now, post-Gödel, the standard interpretation of classical theory seems to be that:

(*f*) PA can, indeed, be interpreted in M;

(*g*) assertions in M are decidable by Tarski's definitions of satisfiability and truth (cf. [Me64], p49-53);

(*h*) Tarskian truth and satisfiability are, however, not effectively verifiable in M[13].

However, the question, implicit in Wittgenstein's argument regarding the possibility of a semantic contradiction in Gödel's reasoning, then arises:

How can we assert that a PA-assertion (whether such an assertion is PA-provable or not) is true under interpretation in M, so long as such truth remains effectively unverifiable in M?

Since the issue is not resolved unambiguously by Gödel in his 1931 paper (nor, apparently, by subsequent standard interpretations of his formal reasoning and conclusions), Wittgenstein's remark can be taken to argue that, although we may validly draw various conclusions from Gödel's formal reasoning and conclusions, the existence of a true or false assertion of M cannot be amongst them.

## 3. What is PA?

### 3.1 Is PA an interpretation of its standard interpretation?

A related philosophical issue is, then, the question:

---

[13] We take this as the interpretation of Tarski's 1936 Theorem (cf. [Me64], Corollary 3.38, p151): The set Tr of Gödel-numbers of wfs of PA which are true in the standard model is not arithmetical, i.e. there is no wf $A(x)$ of PA such that Tr is the set of numbers $k$ for which $A(x)$ is true in the standard model.

Is PA an interpretation of its standard interpretation M?

In other words, since PA is intended to formalise our intuitive arithmetic M as expressed by Dedekind in his formulation of the Peano Postulates, do we accept that PA must be an interpretation of M?[14]

(The standard response to this question seems to lie at the heart of Wittgenstein's reservations, and to be the cause of the uneasiness felt by subsequent philosophers who question the standard interpretations of classical mathematical theory.)

Now, a negative answer would imply that PA cannot be taken as a faithful formalisation of our intuitive, Dedekind, arithmetic; so, either (as standard interpretations of Gödel's reasoning and conclusions implicitly imply) such arithmetic is not formalisable in principle, or there is some interpretation of M that formalises M more faithfully.

The former is, intuitively, an unappealing, Platonic, and implicitly self-limiting, admission; the latter, an unacceptable reflection on the competence of mathematicians to adequately select an appropriate set of primitive, axiomatic, assertions for PA as may be needed for PA to be an effective, and unambiguous, language of precise communication.

An affirmative answer, on the other hand, whilst validating PA as a formalisation of our intuitive, Dedekind arithmetic, would further imply that, since an assertion would then be effectively decidable in PA only if it were effectively decidable in M, there must be some effective method of defining Tarskian satisfiability and truth in M.

---

[14] We note that, prima facie, a strict formalist doctrine would avoid addressing this question; in other words, it would accept PA as a standard formalisation of intuitive arithmetic, but treat both the possible relationship of PA to the parent system that gave it birth, and the significance of the thesis "PA is a standard formalisation of intuitive arithmetic", as being of no mathematical significance for the study of PA as a formal system, and of its interpretations.

## 3.2 Defining effective satisfiability and truth

Although Wittgenstein does not appear to have attempted such a definition - possibly as it may have seemed to involve technicalities beyond the scope of his reflections - we note in §5.1 of Anand [An02c] that such an effective method is, indeed, made available to us by, curiously, a constructive interpretation of Gödel's reasoning and conclusions; an interpretation that is, ironically, more in sympathy with Wittgenstein's constructive approach than Gödel's Platonic one.

## 3.3 Undecidability in PA

Now, the thesis of a constructive interpretation of Gödel's reasoning and conclusions ([An02c], §5.1) is that we may not interpret, for instance, the meta-assertion "PA proves: $[(Ax)F(x)]$[15]" as the non-verifiable, Tarskian meta-assertion:

$F(x)$ is satisfied by any given element $x$ of the domain of M[16].

We must interpret it, instead, as the verifiable meta-assertion:

There is a uniform effective method (algorithm/Turing machine) such that, given any element $x$ of the domain of M, it will effectively decide that $F(x)$ is satisfied in M.

It follows that the meta-assertion, that PA does not prove $[(Ax)F(x)]$, interprets constructively as the meta-assertion:

---

[15] We use square brackets to indicate that the expression within the brackets is to be treated as an abbreviation for an uninterpreted string of a formal system; we use double quotes to indicate that the expression inside the quotes is to be treated as an abbreviation for an interpreted expression in some mathematical language.

[16] Although standard interpretations of classical theory do not appear to highlight this point, this definition implicitly implies that every element of the domain of M (which, by our definition of the language M, would be considered external to M) is necessarily expressible in M.

There is no uniform effective method (algorithm/Turing machine) such that, given any element $x$ of M, it will effectively decide that $F(x)$ is satisfied in M.

Consequently, a constructive interpretation of Gödel's reasoning and conclusions implies that there can be no undecidable assertions in PA; in other words, that PA is syntactically complete!

## 4. Standard interpretations of Gödel's reasoning

### 4.1 How definitive is the standard interpretation of Gödel's reasoning?

However, we are, then, immediately faced with the question: Since the standard interpretation of Gödel's reasoning and conclusions asserts that PA is syntactically incomplete, how definitive *is* the standard interpretation?

Now, in Theorem VI of his 1931 paper, Gödel essentially argues that his "undecidable" proposition, $[(Ax)R(x)]$, is such that (cf. [An02b], §1.3-§1.6):

If $[(Ax)R(x)]$ is P-provable, then $[\sim R(n)]$ is P-provable for some numeral $[n]$.

Now, by standard logical axioms, we have that:

If $[\sim R(n)]$ is P-provable for some numeral $[n]$, then $[\sim(Ax)R(x)]$ is P-provable.

It thus follows that Gödel has essentially argued that:

If $[(Ax)R(x)]$ is P-provable, then $[\sim(Ax)R(x)]$ is P-provable.

Clearly, it should now follow, by the standard Deduction Theorem of first order logic, that:

$[(Ax)R(x) => \sim(Ax)R(x)]$ is P-provable,

and so:

   [~($A x$)$R$($x$)] is P-provable.

## 4.2 We must conclude that PA is omega-inconsistent

However, at this point, standard interpretations of Gödel's reasoning appeal to his explicit assumption that PA is omega-consistent in order to conclude that the PA-provability of [~($A x$)$R$($x$)] cannot be inferred from the above meta-argument.

Now, unless the omega-consistency of PA has some deeper, intuitive significance philosophically, this is not a reasonable inference; since the Deduction Theorem is a fundamental theorem of classical logic, we must, using Occam's razor, conclude from Gödel's reasoning, firstly, that [~($A x$)$R$($x$)] is PA-provable, and, secondly, that PA is omega-inconsistent (and so Gödel's Theorem VI holds vacuously [An02b]).

We note that Wittgenstein's remarks indicate that, prima facie, there appear no intuitively significant philosophical grounds for treating the omega-consistency of PA as a primitive concept; prima facie, there are thus no reasonable grounds for allowing it to over-ride an application of the Deduction Theorem.

In sharp contrast, we note that the omega-inconsistency of PA is natural, and intuitively unobjectionable, under a constructive interpretation of the concept of "PA proves: [($A x$)$F$($x$)]" as described earlier; under such interpretation, an omega-inconsistent PA does not imply that PA, or any of its interpretations, are either inconsistent, or unnaturally consistent (cf. [An02b], Appendix 2). It simply implies that there are arithmetical relations that cannot be verified uniformly by a common algorithm over the domain of their interpretation.

Moreover, as we argue in Anand ([An03c], Appendix 1), this interpretation implies that PA can express relations that are deterministic, yet essentially unpredictable; such a language would have significance for the expression of natural phenomena that are best described in quantum mechanical terms.

The above suggests that it may be the absence of an adequately technical counter-argument that leaves Wittgenstein's viewpoint - and that of others who have shared his reservations - vulnerable to the arguments advanced by standard interpretations of Gödel's reasoning and conclusion; these implicitly imply that any interpretation of Gödel's reasoning and conclusion must be accepted as essentially counter-intuitive on the basis of purely technical considerations.

### 4.3 Rosser's argument may be invalid

Prima facie, the standard interpretation of Gödel's reasoning and conclusions seems strengthened by Rosser's argument that Gödelian undecidability can be established in a simply consistent language without assuming omega-consistency. However, as shown in Anand ([An02a], §7.4(*b*)(*x*)(*3*)), if standard (text-book) expositions ([Me64], p145, Proposition 3.32) of Rosser's proof can be accepted as definitive, then Rosser's argument is invalid.

## 5. When does a formal assertion "mean" what it represents?

An important philosophical issue - which does not seem to have been adequately addressed by standard interpretations of classical theory as yet - is implicit in the key thesis of Floyd and Putnam's paper, and is reflected in Wittgenstein's remark:

> If you assume that the proposition is provable in Russell's system, that means it is true in the Russell sense, and the interpretation "P is not provable" ... has to be given up.

### 5.1 When does a formal assertion "mean" what it represents?

We may state this issue explicitly as:

(*)   When does a formal assertion "mean" what it represents?

Now, if, as argued earlier, we accept that PA formalises our intuitive arithmetic M, and that M is the standard interpretation of PA, it follows that every well-formed formula of PA interprets as a well-defined arithmetical expression, and every well-defined arithmetical expression can be represented as a PA-formula.

The question then arises:

When is an arbitrary number-theoretic function or relation representable in PA?

### 5.2 Formal expressibility and representability

Classically, the question of PA-representability is addressed by the following three definitions (cf. [Me64], p117-118):

(*a*) A number-theoretic relation $R(x_1, ..., x_n)$ is said to be *expressible* in PA if, and only if, there is a well-formed formula $[A(x_1, ..., x_n)]$ of PA with $n$ free variables such that, for any natural numbers $k_1, ..., k_n$:

  (*i*)  if $R(k_1, ..., k_n)$ is true, then PA proves: $[A(k_1, ..., k_n)]$,

  (*ii*) if $R(k_1, ..., k_n)$ is false, then PA proves: $[\sim A(k_1, ..., k_n)]$,

(*b*) A number-theoretic function $f(x_1, ..., x_n)$ is said to be *representable* in PA if, and only if, there is a well-formed formula $[A(x_1, ..., x_n, y)]$ of PA, with the free variables $x_1, ..., x_n, y$, such that, for any natural numbers $k_1, ..., k_n, l$:

  (*i*)  if $f(k_1, ..., k_n) = l$, then PA proves: $[A(k_1, ..., k_n, l)]$,

(*ii*) PA proves: $[(E!l)A(k_1, ..., k_n, l)]$,

(*c*) A number-theoretic function $f(x_1, ..., x_n)$ is said to be *strongly representable* in PA if, and only if, there is a well-formed formula $[A(x_1, ..., x_n, y)]$ of PA, with the free variables $x_1, ..., x_n, y$, such that, for any natural numbers $k_1, ..., k_n, l$:

(*i*) if $f(k_1, ..., k_n) = l$, then PA proves: $[A(k_1, ..., k_n, l)]$,

(*ii*) PA proves: $[(E!l)A(x_1, ..., x_n, y)]$,

### 5.3 When may we assert that $A(x_1, ..., x_n)$ "means" $R(x_1, ..., x_n)$?

We can, thus, re-phrase our earlier question (*) as:

> If a number-theoretic relation $R(x_1, ..., x_n)$ is expressible by a PA-formula $[A(x_1, ..., x_n)]$, when may we assert that, under the standard interpretation, $A(x_1, ..., x_n)$ "means" $R(x_1, ..., x_n)$?

Now we note that, if $R(x_1, ..., x_n)$ is arithmetical, then one of its PA-representation is $[R(x_1, ..., x_n)]$[17], whose standard interpretation is $R(x_1, ..., x_n)$. Hence every arithmetical relation is the standard interpretation of some PA-formula that expresses $R(x_1, ..., x_n)$ in PA, and we can adapt this to give a formal definition of the term "means":

> **Definition 1**: If a number-theoretic relation $R(x_1, ..., x_n)$ is expressible by a PA-formula $[A(x_1, ..., x_n)]$, then we say that, under the standard interpretation, $A(x_1, ..., x_n)$ *means* $R(x_1, ..., x_n)$ if, and only if, $R(x_1, ..., x_n)$ is the standard interpretation of some PA-formula that expresses $R(x_1, ..., x_n)$ in PA.

---

[17] We note that, if $R(x_1, ..., x_n)$ is PA-expressible, then there are denumerable formulas that express it in PA.

The question (*) can now be expressed precisely as:

(**) When is a number-theoretic relation the standard interpretation of some PA-formula that expresses it in PA?

Now, by definition, the number-theoretic relation $R(x_1, ..., x_n)$, and the arithmetic relation $A(x_1, ..., x_n)$, can be effectively shown as equivalent for any given set of natural number values for the free variables contained in them.

However, for $R(x_1, ..., x_n)$ to *mean* $A(x_1, ..., x_n)$, we must have, in addition, that $R(x_1, ..., x_n)$ can be effectively transformed into an arithmetical expression, so that it can be the standard interpretation of some PA-formula that expresses it in PA.

**5.4 What is a mathematical object?**

Clearly, this implies, firstly, that we must be able to add [R] as a primitive, $n$-ary, relation letter to PA, along with suitable axioms, without inviting inconsistency; and, secondly, that $R(x_1, ..., x_n)$ must define a unique mathematical object, where we define (cf. [An02c], §1.2):

(*i*) **Primitive mathematical object**: A primitive mathematical object is any symbol for an individual constant, predicate letter, or a function letter, which is defined as a primitive symbol of a formal mathematical language.

(*ii*) **Formal mathematical object**: A formal mathematical object is any symbol for an individual constant, predicate letter, or a function letter that is either a primitive mathematical object, or that can be introduced through definition into a formal mathematical language without inviting inconsistency.

(*iii*) **Mathematical object**: A mathematical object is any symbol that is either a primitive mathematical object, or a formal mathematical object.

Now, the logical and mathematical antinomies show that we cannot unrestrictedly assume that every arithmetical relation defines a mathematical object[18]. However, if we *can* introduce the *n*-ary relation letter [$R$] into PA as above, without inviting inconsistency, then we can, reasonably, assert that the relation $R(x_1, ..., x_n)$ does, indeed, define a mathematical object in a constructive and intuitionistically unobjectionable way, and that $A(x_1, ..., x_n)$ does mean $R(x_1, ..., x_n)$.

**5.5 Introduction of new symbols by definition**

Now, the question:

> When may we introduce an additional function or relation letter into PA without inviting inconsistency.

is addressed classically by the following proposition (cf. [Me64], Proposition 2.29, p82):

> **Proposition**: Let K be a first-order theory with equality. Assume that K proves: [$(E!u)A(u, y_1, ..., y_n)$]. Let K' be the first-order theory with equality obtained by adding to K a new function letter *f* of *n* arguments, and the proper axiom [$A(f(y_1, ..., y_n), y_1, ..., y_n)$], as well as all instances of the axioms of a first order theory involving *f*. Then there is an effective transformation mapping each well-formed formula *B* of K' into a well-formed formula *B'* of K such that:

---

[18] We note that, by defining a mathematical object precisely, the paradoxical element in the mathematical and logical "antinomies" is effectively eliminated; they define functions or relations that are not mathematical objects. Prima facie, except from a Platonic viewpoint, it thus seems of little significance whether such definitions are taken as defining entities that are mathematically inconsistent (square circle), arguably inconsistent (Pegasus), logically inconsistent (Russell's impredicative set), or mathematical non-objects (the range of Gödel's primitive recursive substitution relation).

(*a*) if *f* does not occur in [*B*], then [*B'*] is [*B*];

(*b*) [(~*B*)'] is [~(*B'*)];

(*c*) [(*B* => *C*)'] is [*B'* => *C'*];

(*d*) [((A*x*)*B*)'] is [(A*x*)*B'*];

(*e*) K' proves: [*B* <=> *B'*];

(*f*) if K' proves: *B*, then K proves: *B'*.

Hence, if *B* does not contain *f*, and K' proves: *B*, then K proves: *B*.

It follows that a number-theoretic function *f* (or relation *R*, if we treat *R* as a Boolean function) may be taken to define a mathematical object if it is strongly representable in PA; we may then introduce the function letter [*f*] into PA without risking inconsistency, and we would then have that $f(x_1, ..., x_n)$ is represented in PA by the formula $[f(x_1, ..., x_n)]$, whose standard interpretation is $f(x_1, ..., x_n)$.

## 6. Why P cannot be interpreted as "P is not provable in PA"

### 6.1 A recursive number-theoretic function that is not a mathematical object

Apart from the trivial resolution of the mathematical and logical paradoxes (*footnote 18*), the significance of the above definitions is that Gödel's primitive recursive substitution function, *Sb*(*x v|y*) (cf. [Go31a], definition 31, p20), is, remarkably, not a mathematical object!

As we note in Anand [An02c]:

> The significance of these definitions is seen in Meta-theorem 1. We prove, there, the existence of an asymmetrical recursive number-theoretic relation - one that is intuitively decidable constructively, but which cannot be introduced through definition as a formal mathematical object into any formal system of Peano Arithmetic without inviting inconsistency; nor, ipso facto, into any Axiomatic Set Theory that models[19] (cf. [Me64], p192-3) such Arithmetic. Hence, it is not a formal mathematical object, and the range of its characteristic function[20] is not a recursively enumerable set[21]!

This is an astonishing result[22]! Firstly, recursive number-theoretic functions and relations are classically accepted as the basic building blocks for defining constructive, and intuitionistically unobjectionable, mathematical objects[23]. Secondly, and in vivid contrast, the relative consistency, and independence, of the Continuum Hypothesis would imply, prima facie, that we may also treat Cantor's non-constructive cardinal, $aleph_1$, as a valid formal mathematical object[24]; thus, we may introduce axiomatic definitions - and an individual constant symbol - for it into any Axiomatic Set Theory without inviting inconsistency!

---

[19] We follow Mendelson's definition of a model ([Me64], p51): An interpretation is said to be a model for a set $T$ of well-formed formulas of P if, and only if, every well-formed formula in $T$ is true for the interpretation in the classical Tarskian sense.

[20] If $R(x)$ is a relation (predicate), then the characteristic function $C_R(x)$ is defined as follows: $C_R(x) = 0$ if $R(x)$ is true, and $C_R(x) = 1$ if $R(x)$ is false ([Me64], p119).

[21] A recursively enumerable set is classically defined as the range of some recursive number-theoretic function, and is implicitly assumed consistent with any Axiomatic Set Theory that is a model for P ([Me64], p250).

[22] Loosely speaking, it may be viewed as a constructive arithmetical parallel to Russell's non-constructive, and seemingly paradoxical, impredicative "set" ([Me64], p2).

[23] See, for instance, Gödel's remarks ([Go31a], p23, footnote 39).

[24] This, essentially, seems to reflect Gödel's point of view, which he expresses in his 1947 paper, "What is Cantor's continuum problem?", whilst discussing whether Cantor's continuum hypothesis should be added to set theory as a new axiom. [Kurt Gödel, Collected Works, vol. 2, Oxford University Press, 1986–2003.]

In other words, Meta-theorem 1 in Anand [An02c] establishes that there is no PA-formula that can mean $Sb(x\ v|y)$ under interpretation.

## 6.2 No PA-formula P can interpret as "P is not provable in PA"

Now we note that Gödel's number-theoretic relation $Bew(x)$ - which holds if, and only if, $x$ is the Gödel-number of some PA-provable formula P - is defined ([Go31a], definition 46, p22) in terms of his primitive recursive number-theoretical relation $Sb(x\ v|y)$. It follows that, since the latter is not a mathematical object, neither is the former.

Thus, there is no PA-formula that can mean $Bew(x)$ under interpretation; ipso facto, there is no PA-formula that can interpret as an arithmetic assertion that is equivalent to the assertion "P is not provable in PA". As we remark in Anand ([An02c], §II(3)), such an assertion can only be a convention, not an arguable inference.

## 6.3 No PA-formula can interpret as "PA is consistent"

It also follows that if the concept of "PA-provability", as defined by Gödel, cannot be expressed in PA, then Gödel's Theorem XI [Go31a], which asserts that the consistency of PA cannot be established within a consistent PA, must hold vacuously. This argument is discussed in detail in Anand ([An02c], §II), where we prove Theorem 2, to the effect that no PA-formula can mean "PA is consistent" under the standard interpretation.

We essentially argue there that:

> If we assume that the number-theoretic sentence $(Ex)(Form(x)\ \&\ \sim Bew(x))$ ([Go31a], p36, footnote 63), abbreviated $Wid(PA)$, defines the proposition "PA is consistent" in a constructive and intuitionistically unobjectionable way[25], then can we consistently assume further that *Wid*(PA) is equivalent to the standard interpretation of some PA-formula [*Con*(PA)]?

---

[25] In other words, this definition can be assumed equivalent to Mendelson's classical meta-definition of consistency ([Me64], p37). We argue in Anand ([An02c], §II-2), albeit in a different context, that this assumption, too, may conceal an implicit ambiguity.

... the latter assumption is one of the implicit meta-theses that appear to underlie Gödel's proof of, and the conclusions he draws from, his Theorem XI ([Go31a], p36).

Clearly, the reasoning in Meta-theorem 1 and Meta-lemma 2 implies that such an assumption is invalid[26]. We conclude that there is no PA-formula, *Con*(PA), whose standard interpretation is the number-theoretic assertion (E$x$)(*Form*($x$) & ~*Bew*($x$)) - which is defined by Gödel as equivalent to "PA is consistent".

## 7. Conclusion

We conclude that standard interpretations of Gödel's reasoning and conclusions are not definitive. Wittgenstein's reservations on Gödel's interpretation of his own formal reasoning, as reflected in his "notorious" paragraph, and the uneasiness that academicians and philosophers such as Floyd and Putnam continue to sense, and express, over standard interpretations of Gödel's formal reasoning, must be respected as significant indicators of possible ambiguities that may be rooted in implicit assumptions underlying standard interpretations of classical foundational concepts.

---

[26] For, if (E$x$)(*Form*($x$) & ~*Bew*($x$)) is the standard interpretation of one of its formal representations, then so also are each of *Form*($x$) and *Bew*($x$). Arguing similarly, this would eventually imply that the recursive function *Sb*($x$ 19|*Z*($y$)) too is the standard interpretation of one of its formal representations, contradicting Meta-theorem 1.

*Acknowledgement: The idea of relating the uneasiness expressed by Wittgenstein, Floyd, Putnam and others - over standard interpretations of classical theory - to a constructive interpretation was conceived after viewing the recent discussion in the Foundations of Mathematics Forum, FOM, over Floyd and Putnam's paper [FP00].*

*Author's e-mail: anandb@vsnl.com*

(*Updated: Thursday 29$^{th}$ May 2003 6:42:51 PM IST by re@alixcomsi.com*)